\documentclass[table, svgnames, dvipsnames, letterpaper, 11pt]{report}

\usepackage{amsmath, amsthm, latexsym, amssymb, graphicx, graphics, bold-extra, multirow,
pdfpages, array, float, subcaption, bm, changepage, longtable, booktabs, algorithmic, wrapfig}

\usepackage{xcolor, soul}

\usepackage{mathrsfs}
\usepackage{stmaryrd}

\usepackage{titlesec}
\titlespacing\section{0pt}{12pt plus 4pt minus 2pt}{0pt plus 2pt minus 2pt}
\titlespacing\subsection{0pt}{12pt plus 4pt minus 2pt}{0pt plus 2pt minus 2pt}
\titlespacing\subsubsection{0pt}{12pt plus 4pt minus 2pt}{0pt plus 2pt minus 2pt}

\usepackage{etoolbox} 
\patchcmd{\endabstract}{\null}{}{}{}
\patchcmd{\thebibliography}{\chapter*}{\section*}{}{}
\makeatletter
\patchcmd{\@tocline}{\hfill}{%
  \leaders\hbox{$\m@th
    \mkern \@dotsep mu\hbox{.}\mkern \@dotsep
  mu$}\hfill}{}{}
\makeatother

\makeatletter
\patchcmd{\l@section}{\hfil}{%
  \leaders\hbox{$\m@th
    \mkern \@dotsep mu\hbox{\color{blue}.}\mkern \@dotsep
  mu$}\hfill}{}{}
\makeatother


\usepackage{tocloft}

\usepackage[T1]{fontenc}
\usepackage{lmodern}
\usepackage{textcomp}

\usepackage{arydshln}

\usepackage[sort&compress, square, numbers]{natbib}

\usepackage[resetlabels]{multibib}

\usepackage[title]{appendix}
\usepackage[export]{adjustbox}

\usepackage{tikz} 
\usetikzlibrary{positioning}

\usepackage{geometry}
\theoremstyle{remark}

\usepackage{fancyhdr}
\usepackage[pdftex, breaklinks=true]{hyperref}
\hypersetup{
	unicode=false,								
	pdftoolbar=true,							
	pdfmenubar=true,							
	pdffitwindow=true,							
	pdftitle={ACMS Jorunal and Proceedings},	
	pdfauthor={Russell W.\ Howell},				
	pdfsubject={Subject},						
	pdfnewwindow=true,							
	pdfkeywords={keywords},						
	colorlinks=true,							
	linkcolor=blue,								
	citecolor=blue,								
	filecolor=magenta,							
	urlcolor=blue,								
	linktoc=all
}

\makeatletter
\renewenvironment{abstract}{%
      \@beginparpenalty\@lowpenalty
      \small
      \begin{center}%
        \bfseries \abstractname
        \@endparpenalty\@M
		\vspace*{-0.2in}
      \end{center}\quotation}%
     {\endquotation\par%
     }
\makeatother

\begin{document}

\fancyhf{} 
\lfoot{\bf \footnotesize ACMS \emph{Journal and Proceedings}}
\cfoot{\bf -- \thepage~--}
\rfoot{\bf \footnotesize 23rd Biennial Conference}

\setcounter{section}{0}
\setcounter{subsection}{0}
\setcounter{subsubsection}{0}
\setcounter{footnote}{0}
\setcounter{figure}{0}
\setcounter{lemma}{0}
\setcounter{equation}{0}
\setcounter{theorem}{0}
\setcounter{table}{0}

\label{mangum}

\begin{center}
{\bf \Large The Nature of Reasoning in \\ Theology, Philosophy, and Mathematics}

\href{https://www.clemson.edu/science/academics/departments/mathstat/about/profiles/crmangu}{Chad R.\ Mangum}\,  
\href{https://www.clemson.edu/}{(Clemson University)}

\fbox{
\begin{minipage}[t]{1.03in}
\vspace*{0pt}
\includegraphics[height=1.26in]{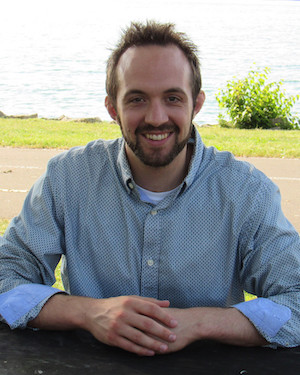}
\end{minipage}
\begin{minipage}[t]{4.5in}
\vspace*{0pt}
Chad R. Mangum (Ph.\ D., North Carolina State University) teaches mathematics at Clemson University. His
research interests include Lie algebra representation theory and related areas, as well as apologetics,
philosophy, theology, music, and graph theory. He deems novel connections between his areas of interest
to be especially intriguing.  Chad also enjoys drumming, playing sports, and staying involved in his
church.
\end{minipage}
}
\end{center}

\begin{abstract}
\noindent
This article supports the epistemological claim that sound human reasoning about ultimate
knowledge is either foundational or circularly justified. In particular, questions which naturally arise in theology,
philosophy, and related disciplines, to the extent that they rationally treat ultimate knowledge, are
necessarily supported in one of these ways. Comparisons with, contrasts to, and analogies from
mathematics are given to illustrate and enhance this central claim.
\end{abstract}

\section{Introduction}\label{introduction}
One of the most basic goals of any rational human enterprise is to intellectually justify one's truth
claims. This article will investigate the nature of that justification process. We have in mind what we will
call ``propositional justification,'' that is, justification \emph{of} propositions \emph{by} propositions 
(truth-valued statements). 
Words like ``support'' or ``justification'' will refer to this notion of propositional justification unless 
specified otherwise. 

Our aim is to examine the \emph{structure} of rational justification, not to arbitrate rationality itself. 
The question we seek to answer is not ``is proposition $P$ rationally held?'' but rather, ``\emph{given} a theory 
of warrant,\footnote{A set of (and relations between) standards or criteria by which rationality is measured; it goes by other names as well.} what is the 
nature of the structure of the rational support of proposition $P$?''

It is natural to relate these epistemological issues to mathematics as the latter is commonly regarded
as occupying rarefied air, perhaps even uniquely, in its theoretical completeness. Famously adored as the
``queen of the sciences''\footnote{While applying this term to mathematics is attributed to Carl
Friedrich Gauss, math was not its original object. It referred to theology (of all disciplines!) during
the Middle Ages. See \cite[p.~213]{BH}.} with the ``most uncontroversial examples of proof,''
\cite[p.~56]{FR} if something is observed in mathematical thinking, it is often expected to apply to
other (maybe even all) rational endeavors. Examining the interplay between mathematics and ultimate
justification for truth claims is therefore worthwhile and will be taken up in this article.

This article is split into three sections. Section \ref{theargument} states and argues for the main thesis.
Section \ref{connectmath} uses mathematics as a case study to make connections between the notions of foundational knowledge and 
circular reasoning, particularly noting how circularity is avoided by the
major schools of thought in the philosophy of math.  Section \ref{followup} briefly addresses some
questions which naturally follow.

\section{The Argument}
\label{theargument}

\subsection{Assumptions and Definitions}
It would be a Herculean (or impossible) task to fairly and incontrovertibly list \emph{every} assumption;
complete, untainted epistemological self-awareness might be only a theoretical ideal. We will nonetheless
attempt to identify the ones most fundamental and most relevant to the present topic.

Our primary object in this section will be the ``hypothetical belief set'' ($HBS$) of a fixed human agent $H$ at a fixed
time $t$. It is ``hypothetical'' because we allow the set to include beliefs that $H$ \emph{would}
hold if prompted appropriately, but which $H$ might not otherwise believe or even conceive of. The set
$HBS$, sometimes written $HBS_t$ to emphasize the dependence on $t$, along with the relations among
its elements, is the \emph{noetic structure} of $H$. \cite[p.~48]{PL} We will assume $H$'s thoughts are not 
past eternal and that it takes some positive amount of time for $H$ to form a thought. Coupled 
with the dependence on $t$, we may thus assume that the cardinality of $HBS_t$ is finite.

We will assume the laws of logic and general principles of rationality which are too many (and, likely,
too familiar) to enumerate. But as they will be essential terms, let us define them a bit more
thoroughly.\footnote{We take the notions of rationality and logic to be \emph{rationally} foundational
(as seems to be required of rational human endeavors), but not ``semantically foundational.'' 
As will be discussed throughout the article, there are two sets in any rational human endeavor
to which our thesis applies: rational support and definitions. An element which appears in both sets need
not be ultimate (in the sense defined later in the section) in both sets, modulo the caveat discussed
below that definitions can be converted to truth claims. In this article, ``rationality'' is considered
ultimate in the former set, but not the latter. Nonetheless, the definitions we give should be taken to
be informal.} Here, we will consider logic to be formal, classical logic.\footnote{As opposed to the
informal use of the word ``logic,'' which can refer to any kind of reasoning, and as opposed to
nonstandard formal logics. If one insists on a nonstandard logic, it is likely to be unproblematic for
this article.} Definitions include

\vspace*{-0.2in}

\begin{itemize}
\item
  ``a formal, scientific method of examining or thinking about
  ideas;''\footnote{\url{https://dictionary.cambridge.org/us/dictionary/english/logic},
    accessed June 20, 2022.}
\item
  ``a science that deals with the principles and criteria of validity of
  inference and demonstration\textbf{:} the science of the formal
  principles of reasoning;''\footnote{\url{https://www.merriam-webster.com/dictionary/logic},
    accessed June 20, 2022.}
\item
  ``the science of thinking about or explaining the reason for something
  using formal methods.''\footnote{\url{https://www.oxfordlearnersdictionaries.com/us/definition/english/logic_1?q=logic},
    accessed June 20, 2022.}
\end{itemize}

\vspace*{-0.2in}

We consider deduction to be the process of drawing conclusions by applying the rules of logic to 
one or more statements (elements of $HBS$). Rationality is

\vspace*{-0.2in}

\begin{itemize}
\item
  ``the quality of being based on clear thought and reason, or of making
  decisions based on clear thought and reason;''\footnote{\url{https://dictionary.cambridge.org/us/dictionary/english/rationality},
    accessed June 20, 2022.}
\item
  ``the quality or state of being agreeable to reason;''\footnote{\url{https://www.merriam-webster.com/dictionary/rationality},
    accessed June 20, 2022.}
\item
  ``the fact of being based on reason rather than emotions.''\footnote{\url{https://www.oxfordlearnersdictionaries.com/us/definition/english/rationality},
    accessed June 20, 2022.}
\end{itemize}

\vspace{-3mm}
Accordingly, we will call an intellectual endeavor or claim \emph{rational} if it is ultimately supported
by reasons. What constitutes a ``reason'' is a criterion which we do not attempt to define in this
article (beyond noting that a reason must be a truth-valued statement). That is, we will assume that some
theory of warrant has already been accepted; there is much epistemological
debate on such topics on which we (attempt to) remain mostly silent. We simply consider whether a claim
is supported by reasons or not, without regard to whether they are ``good reasons.'' Note also that there
is a contrast with rational\emph{ism}, the philosophical view that ``regards reason as the chief source
and test of knowledge.''\footnote{See \url{https://www.britannica.com/topic/rationalism}, accessed June
20, 2022.}

We consider, then, \emph{irrationality} to be the negation of rationality and \emph{irrational} the
negation of rational; in particular, for brevity, the notion of ``arationality'' (being outside the arena
of reason, even if not in conflict with it) is subsumed under an umbrella term of irrationality in this
article. Academic disciplines such as theology (the study of religious
belief), religion (the pursuit of the supernatural or supreme existence), and philosophy (the study of
the foundational nature of reality) will be referred to in broad strokes.

By \emph{circular reasoning}, therefore, we will mean the support of a proposition in $HBS_t$ for which the conclusion is
among the premises. We will use the term ``circularity'' for brevity to refer to circular human reasoning. 
It is sometimes
referred to in the literature as ``begging the question'' or ``petitio principii,'' though some authors
make a distinction among these terms. 
A proposition will be considered \emph{foundational} if it has no propositional support in $HBS_t$.

We will assume an agreed-upon meaning\footnote{Or ``approximate meaning.''} of all terms used, with only
a few key terms being explicitly defined. Of course, this is a nearly universal assumption, and any
meaning-carrying dialogue depends on it, so it is not overly ambitious. Nonetheless, given the
epistemologically fundamental nature of the content that follows, adequate epistemological self-awareness
dictates that this assumption be stated, along with the acknowledgement that semantic misunderstandings
are still possible regardless. It is the author's hope in defining a few key terms that any such
persisting misunderstandings are insignificant to the central thesis.

We fix a claim, $P \in HBS_t$, of genuine knowledge about reality, i.e. $P$ is (agreed upon to be) true.
Likewise, we assume a goal of rational investigation of objective truth about reality, as opposed to
statements grounded ultimately in opinion, preference, feeling, or other subjective standards. In
particular, issues raised by the philosophical positions of skepticism, which doubt whether truth can be
known or if there is any objective reality at all, will not be considered herein. Thus, we will assume
that $H$ can justifiably claim to know \emph{something}. Note also that definitions could be rephrased
as truth claims in the form ``term $X$ has meaning $Y$;''\footnote{A similar phenomenon can be seen in
mathematics, when sometimes definitions are referred to as axioms; indeed, ``there is little distinction
between a mathematical axiom and a definition.'' \cite[p.~216]{BH}} $P$ might be as such.

It is useful to refine the epistemic \emph{level} of $P$: either $P$ is ultimate knowledge
or it is not. In this article, \emph{ultimate knowledge} is the most basic knowledge that $H$ possesses. 
We will often refer to such a proposition as \emph{an} ultimate, or ultimates in plural (using the
word ``ultimate'' as a noun). An ultimate is something for which $H$ does not have a
deeper knowledge \emph{level} to undergird it.  Likewise, we will call questions about ultimate knowledge ``ultimate
questions.'' What we call ``ultimates'' in this article are given many other names depending on the context and
author: presupposition, basic commitment, \cite[p.~xvii]{FR} (properly) basic belief, \cite{PL},
\cite[p.~59]{FR}, starting point, innate belief, final authority, foundation, foundational belief,
assumption, and so on (and, unfortunately, sometimes these terms are used with different meanings).\footnote{At $t$, we consider the set of foundational propositions to be a subset of the set of ultimates.}  It
is tempting to equate this notion with the mathematical notion of axiom (or ``hypothesis'' to some), but
as Christian theologian John Frame discusses, \cite[pp.~xxxi-xxxii]{FR} and as we shall see shortly,
there is a distinction.

The aforementioned assumptions are not necessarily themselves ultimates, but are a practical
necessity to focus the discussion and keep it to a manageable length.

\subsection{Primary Thesis: Ultimate Justification for a Truth Claim}

What we discuss below is a version of Agrippa's trilemma, credited to the philosophical skeptic Agrippa
and the Pyrrhonists.\footnote{The essentials of the problem go by other names as well, such as
Münchhausen's trilemma, or Fries' trilemma. The term ``Münchhausen's trilemma'' is credited to \cite{AL}.
The present article uses the framework of this trilemma, while coming to different conclusions than
Agrippa and the Pyrrhonists.} It is one manifestation of the problem in epistemology regarding ultimate
justification for a truth claim. The nuanced differences of its disparate treatments will not be
essential here.\footnote{For example, in a work by Greek philosopher Sextus Empiricus, five possible
outcomes to the question of ultimate justification are given. Three of those five are contained in
Agrippa's trilemma, with two additional ones which are disallowed by our current assumptions:
dissent (in our context, this would be to assume that the proposition $P$ is not agreed to
be true) and relation ($P$ appears to be true from the point of view of a certain individual,
but that does not necessarily constitute a statement of the inherent truth of $P$; in other words, it is
the conclusion that truth does not refer to an objective reality).
See \cite{EM}. For a readable overview of the issues in Agrippa's trilemma, 
refer to the Stanford Encyclopedia of Philosophy's entry on Skepticism; see
\url{https://plato.stanford.edu/entries/skepticism/} (accessed July 31, 2022).}

Our primary thesis is the claim that 
\begin{center}
(FC) $\hspace{3mm} HBS_t$ contains elements which are foundational or are circularly supported.
\end{center}
We will begin with $P$ to conclude (FC). One can ask for proximate, rational justification of $P$: 
why, or on what rational authority (empirical
evidence, logic, the authority of the source from which $P$ was discovered, or whatever reason might
be legitimate according to $H$'s theory of warrant), can we say that $P$ is the case?\footnote{We can
consider such a question as a ``Chisholm-like question'' without the misunderstanding
of how such a question should be taken as in \cite[p.~51]{PL}.}
Suppose a distinct supporting reason, $Q$, is given.
One could equally ask for the distinct, proximate, rational justification of $Q$, call it $R$. The pattern
continues, descending into a sequence of ``why'' questions and answers: $P, Q, R, \ldots$. We will call
this the \emph{support sequence}.\footnote{Could $P$ have a distinct support sequence, say $P, E, F, G, \ldots$?
We allow for that possibility. The analysis which follows will fix a particular sequence,
though the conclusions apply to any support sequence which one desires to choose. The salient
characteristic is that $P$ has only finitely many distinct support sequences, which follows since 
$HBS_t$ is assumed to be finite in size. Thus, fixing just one support sequence suffices.} Under our
assumptions, there are only two options\footnote{The reader may wonder where the third horn of Agrippa's
trilemma is. Two of them are together contained in a single option: the finite termination of the
sequence.} with regard to the termination of this sequence (and hence, for the foundations of 
$H$'s noetic structure) which we will explore in turn:

\begin{enumerate}

\vspace*{-0.2in}

\def\labelenumi{\arabic{enumi}.}
\item
the sequence never terminates, continuing to give distinct supporting reasons; or

\item
the sequence terminates in finitely many steps, ceasing to give distinct supporting reasons.
\end{enumerate}

\vspace*{-0.2in}

\subsection{Option 1: The ``Endless Why''}
This option implies there is an infinite regress of supporting reasons for claim $P$, a truly endless
sequence of ``why'' questions and answers.  In Agrippa's trilemma, this is the ``infinitism'' horn. 
This contradicts our assumption that $HBS$ is finite at $t$, and therefore can be discarded.

Additionally, perhaps enlightening is the observation that, even if the finiteness assumption were dropped,
it is questionable if $P$ could be considered to be genuinely rationally supported. If the argument
supporting $P$ is never complete, is $P$ really supported after all? Is an in-progress proof really
a proof (in any theory of warrant)? One particularly troublesome potentiality: the support sequence could 
eventually contain something
that breaks our assumptions, and, even if that hasn't occurred yet, we may not have sufficient
confidence that it won't occur at one of the future steps. Such observations demonstrate that, even under
slightly more general assumptions, an infinite
support sequence for $P$ is dubious at best, and possibly irrationality ``in disguise.'' 

\subsection{Option 2: Foundational or Circular Support}
Therefore, the support sequence cannot be truly endless in rational support of $H$'s genuine knowledge;
it must terminate in finitely many steps. That is, some supporting reason $Z$\footnote{The termination of the
support sequence of $P$ is likely to be a set of reasons rather than one single statement. For ease of notation,
and since the individual elements of such a set need not be separated for our purposes, we will simply notate
the termination $Z$ and treat it as a single reason.} must be the last distinct supporting reason at $t$ (where we allow
for the sequence to terminate in 0 steps, \emph{i.e.}, $P = Z$). That $Z$, by definition, 
has no deeper or more fundamental level with distinct reason(s) supporting it at $t$, so it is at the ultimate justificatory level.
How does $Z$ relate to other elements of $HBS_t$?

That question is answered in two different ways in Agrippa's trilemma, the first being
\emph{foundationalism}.\footnote{The related stance of positism, \cite{EG} in which $Z$ is simply unjustified
(presumably relative to the chosen theory of warrant), is discarded as it contradicts our assumptions.
In being unjustified by its own admission, it represents a choice that is not rational truth-seeking.} 
Here, $Z$ is accepted as foundational without further \emph{propositional} rational defense 
(though it may be defended some other way). The support sequence merely stops at $Z$.
Depending on the nature of $Z$ and the theory of warrant,
$Z$ could be considered knowledge without being inferred from \emph{any} proposition, including from itself. \cite{PL}

The other horn of the trilemma is referred to as \emph{coherentism} in which $Z$ is part of a mutually 
supporting, self-contained theory. This requires a ``cycle'' of reasons;\footnote{Note that $Z$ is the final
\emph{distinct} reason. It could be the case in $HBS_t$ that the reason given for $Z$ is some other proposition
$R$ already listed in the support sequence, thereby initiating cyclical traversal through the sequence.}
this is circularity. 

Therefore, there are only two non-rejected options for the ultimate justification of $P$:\footnote{That is, only
two options which do not contradict our assumptions.} foundational support and circular support.\footnote{The
interplay between them is examined deeply in \cite{HA}.}
Thus we conclude (FC).

Most connections to mathematics will be saved for Section \ref{connectmath}, but one is appropriate here.
We may impose a finite directed graph structure on $HBS_t$ since we view it as a finite, discrete set 
with nonsymmetric relations among (some) pairs of elements. Call this the ``knowledge graph'' of $H$ and label it $G$.
\footnote{This may be a metaphor to add to the list of ``raft'' and ``pyramid.'' See \cite{SO}.}
Nodes of $G$ are propositions (elements of $HBS$), and arrows point from one node to the others which support it. Every
node in a connected\footnote{Or, more commonly, ``weakly connected.'' See \cite[p.~163]{CZ}.} component
of $G$ has at least one arrow pointing to it or emanating from it. A
theorem\footnote{\cite[p.~89]{CZ}, Theorem 4.3 proves this result for an undirected graph. The directed
version can be found in Lemma 1 at \url{https://www.math.cmu.edu/~af1p/Teaching/GT/CH10.pdf}, accessed
June 15, 2022.} in graph theory tells us that every such graph contains a cycle
(which would indicate circular reasoning in our analogy) or a sink (a node with no arrows emanating from
it; in our analogy, this is a foundational proposition).\footnote{Note, then, that a foundational proposition cannot be part of a cycle, but an ultimate might be as per (FC).} If one accepts the knowledge graph analogy, then
the proof of this theorem offers an alternative proof of (FC). Perhaps this mathematical analogy can be mined to 
provide further insight into epistemology.

%

\subsection{Why entertain the plausibility of circular reasoning at all?}\label{entertaincirc}
Nearly all rational thinkers eschew any form of circular reasoning as logically fallacious.
Philosopher Douglas N. Walton, with an entire book devoted to circular reasoning, \cite{W1} explains
succinctly why this is the case: ``Arguing in a circle becomes a fallacy by basing it on prior acceptance
of the conclusion to be proved. So the fallacy of begging the question is a systematic tactic to evade
fulfillment of a legitimate burden of proof.'' \cite[p.~66]{W2}, \cite[pp.~254-5]{FR} As philosopher S.
Morris Engel puts it, ``if the supporting premises merely repeat or rephrase what is stated in the
conclusion, as in all cases of begging the question, the argument contains no premises and is therefore
fallacious.'' \cite[p.~147]{EN}, \cite[p.~255]{FR} According to theologian Joseph E. Torres, ``{[}i{]}f
question-begging is embraced, fideism is implied. Fideism is the rejection of a rational
{[}argument{]}.'' \cite{TO}, \cite[p.~255]{FR}\footnote{Fideism, beliefs taken on faith alone without any
rational support, is widely seen as anti-intellectual.}
Even in everyday discourse, labeling a putatively intellectual argument to be circular is a death blow.

Why, then, allow for circular reasoning as a plausible option at all?
The short answer is that the rationality of circular (or any other kind of) reasoning is left
up to the theory of warrant, and in this article, we attempt not to evaluate such theories.
But there is another layer worth examining.

A variety of circularity is discussed, and evidently promoted, by Christian presuppositional apologists, 
a school of thought expounding
on the works of theologian Cornelius Van Til. For example, John Frame defines ``circularity'' as ``an
argument in which the conclusion justifies\footnote{It is not entirely clear if Frame's notion of ``justification'' is the \emph{propositional} variety. Similarly, his notion of ``ultimate'' seems to refer to a standard or criterion in the theory of warrant rather than to a proposition.} itself. All arguments seeking to prove the existence of an
ultimate or final authority are circular in this sense.'' \cite[p.~291]{FR} Philosopher and apologist Greg Bahnsen, often considered
the primary popularizer of Van Til's work,\footnote{Rev. Dr. K. Scott Oliphint, 
expresses this in his
comment of \cite{B1}. Oliphint says, ``For those who want to understand Van Til, whether to agree or
disagree, at least two things are both essential and too often neglected. The first is to read Van Til,
the second is to read Greg Bahnsen.'' (Back cover).} says that this kind of circularity is ``involved in
a coherent theory (where all the parts are consistent with or assume each other) and which is required
when one reasons about a precondition for reasoning.'' \cite[p.~518n122]{B2}
Theologian R. C. Sproul concurs: ``That all reasoning is ultimately circular in the sense that
conclusions are inseparably related to presuppositions is not in dispute.'' \cite[p.~70]{SP},
\cite[p.~258]{FR} Seeing circularity as a consequence of the finiteness of humanity is expressed
succinctly by Van Til himself, who writes, ``{[}w{]}e hold it to be true that circular reasoning is the
only reasoning that is possible to \emph{finite} man'' (emphasis added). \cite[p.~12]{V2},
\cite[p.~254]{FR}

How can serious thinkers seem to accept an argumentation tactic that is clearly fallacious?
It seems that they refer to different
types of rational support. When seeking \emph{propositional support}, as in this article, there is broad agreement,
even among presuppositionalists, that circular reasoning is fallacious. The variety of circularity endorsed 
by presuppositionalists is of a different type, a type we will call \emph{authoritative support},\footnote{Perhaps the term ``foundational support'' would be equally descriptive.} referring to the
\emph{criteria} which one uses to judge veracity or rationality. Frame recognizes the necessity:
``{[}W{]}hen one is arguing for an ultimate criterion, \ldots{} one must use
criteria compatible with that conclusion. If that is circularity, then everybody is guilty of
circularity.'' \cite[p.~11]{FR}\footnote{Indeed, ``reasoning in a vicious circle is the only alternative to reasoning in a
circle.'' \cite[p.~24]{V1}, \cite[p.~257]{FR} In the words of philosopher Richard Pratt, ``in
argumentation, reasoning should be linear,'' but circularity of the ``\emph{process} by which finite minds
attain knowledge to be used in arguments'' is ``inescapable.'' \cite{PR}, \cite[p.~257]{FR} (emphasis added).} In other words, rhetorically, by what standard can an \emph{ultimate criterion} of
knowledge be aptly declared suitable, if not by itself (for, if not by itself, then the criterion is not ultimate)?\footnote{It is in this sense that ultimates occupy a ``privileged position'' in rational argumentation in that they exhibit
special epistemological properties not seen in non-ultimates. See Bahnsen's lecture ``Reasoning With Unbelievers,'' beginning at 25:20:
\url{https://www.youtube.com/watch?v=6RBz-zAWoKk}, accessed December 8, 2021. See also the discussion of the ``meta'' relationship that ultimates have to the support sequence in Sections \ref{lackECtheophil}, \ref{circinarticle}.}
By what criteria of knowledge shall we judge our criteria of knowledge but themselves?
To emphasize the distinction between fallacious circular reasoning of propositional support and purported circularity of authoritative support, Torres terms the latter
``circular coherence.'' \cite{TO}, \cite[p.~256]{FR}\footnote{More from Torres on this (and related
topics) can be found at ``Presuppositionalism and Circularity\ldots Again?'' at
\url{https://apolojet.wordpress.com/2012/03/15/presuppositionalism-and-circularity-again/}, accessed June
20, 2022. A similar bifurcation of types of ultimates, propositional versus authoritative, might be in order as well.}
This is the version of circularity endorsed by the presuppositionalist school.

It is worth emphasizing that this type of circularity and the question begging fallacy (petitio principii) are distinct notions. According to Torres, a presuppositionalist, ``Van Tillians, at least implicitly, distinguish between circular coherence and begging the question, embracing the former and rejecting the latter.'' \cite[p.~258]{FR}\footnote{Other terms are used for the distinction
between circular coherence and question begging: Frame calls the former circularity
``broad'' and the latter ``narrow.'' \cite[p.~14]{FR} Torres labels them ``virtuous'' and ``vicious,'' respectively. \cite{TO}, \cite[p.~255]{FR} 
As Torres explains, ``the main distinguishing mark between these
two forms of circularity is how one handles the burden of proof, or evidential priority.'' \cite{TO},
\cite[p.~264]{FR} Thus, the question-begging fallacy is ``vicious,'' whereas circular coherence is seen as ``virtuous.'' (cf. Russell's ``vicious circle principle'' discussed in Section
\ref{connectmath}.)} The circularity of presuppositionalists is a requirement of the theory of warrant to use only reasons within argumentation (arrows in $G$) that comport with the ultimates (each would appear as a node in $G$, either a sink or as part of a cycle); thus it is some relation between arrows and ultimates (specific nodes) which does not appear explicitly in $G$ (in particular, presuppositional circularity does not require a cycle in $G$). A type of ``circle'' appears in that ultimates and the theory of warrant are self-consistent,\footnote{That is, the ultimate criteria underlie the reasoning used, even in (linear) argumentation for their own legitimacy, in a way that non-ultimates do not. This ``circle,'' however, may be non-propositional in nature, and hence not be circular reasoning as traditionally understood.} but this is \emph{not} a cycle in $G$.
Question begging is a relation between various nodes of $G$ and \emph{is} evident by the presence of a cycle in $G$. As the presuppositional notion of circularity is distinct from the familiar concept of circular reasoning, criticism of presuppositionalists for ambiguity in language may well be warranted.\footnote{Frame recognizes the risk: ``I don't care very much whether the Christian apologist accepts or rejects the term \emph{circular} to describe his argument. There are obvious dangers of misunderstanding in using it...I am more inclined now to say to my critics, `Granted your definition of circularity, I don't believe in it.' " \cite[p.~11n19]{FR}}

 

Thus the contrast of support type (propositional versus authoritative) demonstrates how the presuppositional
notion of circularity is distinct from and not contradictory to the common understanding which rightly disparages circular reasoning.
As authoritative support reveals itself in evaluating rationality and not the structure among propositions,
questions of authoritative support reside in the theory of warrant.
Hence, we do not judge such positions here, but underscore the point that permissibility of circularity in \emph{authoritative} support is not inconsistent with fallaciousness of circularity in \emph{propositional} support.



\section{Connections to Mathematics}\label{connectmath}
\subsection{The Epistemic Characteristics}\label{TheECs}

It is justifiable to bring mathematics into this discussion so long as we assume that at least one
statement of mathematics is genuine knowledge (to play the role of $P$ from Section \ref{theargument}). We
add that assumption for this section.

Mainstream mathematics today follows the axiomatic method.\footnote{The discipline settled on this more
structured approach ``since the free use of our conceptions\ldots led to disaster,'' specifically the
paradoxes that arose from Cantorian set theory, discussed below.  \cite[p.~40]{KE}} It
begins\footnote{That is, the theoretical construct begins in the way described. The process in which
humans engage -- the process of doing mathematics -- rarely proceeds in this order chronologically.} with
a set of (formally undefined) terms\footnote{This idea has a long heritage in mathematics.  As far back
as Aristotle it was acknowledged that some terms ``must be undefined or else there would be no starting
point.'' \cite[p.~20]{KI} (Note the lurking foundationalism!) Interestingly, Euclid's \emph{Elements} infamously attempted to define all
concepts used, despite there being historical indications that he was aware of Aristotle's works.
Subsequently, the work of Euclid had an impact on ``practically all the mathematicians who followed him
for two thousand years.'' \cite[pp.~101-2]{KI} We will comment further on this in the subsection on
formalism.} and axioms (asserted-but-unquestioned assumptions, taken to be true without
proof),\footnote{Formally, one could view an axiom as proving itself.} and uses logical deduction, enhanced by
intuition,\footnote{One might prefer the word ``creativity'' instead, or any of several related ideas.}
to justify new statements. For convenience, we will give names to two epistemic characteristics exhibited
by this process of mathematics. Mathematics

\vspace*{-0.2in}

\begin{enumerate}
\def\labelenumi{\arabic{enumi}.}
\item
(EC1: starts somewhere) begins with foundational principles, recognizing that mathematical knowledge must
start with something;\footnote{This is expressed nicely by $18^{th}$ century Scottish philosopher Thomas Reid
in \cite{CL}, esp. p. 148.} and

\item
(EC2: restricts the scope) restricts what is properly a ``mathematical question'' and what is not. Some
questions, even some which refer to indisputably mathematical objects, are not mathematical questions. In
particular, there are logical arguments which are not part of mathematics (and are usually placed within
the realm of theology and/or philosophy) because those arguments refer to nonmathematical questions or objects.
\end{enumerate}

The content of mathematics is whatever can be concluded (according to the methods allowable in the given
philosophical school, as discussed below) from the starting point(s) of EC1 and within the bounds of EC2.
The process of mathematics can be seen through the lens of propositional support as discussed
in Section \ref{theargument}.
On the basis of the elements of EC1 (principally axioms and definitions), mathematics is frequently 
seen as being foundationalist, \cite[p.~215]{BH}, \cite[pp.~51-2]{HB} even if artificially so vis-\`{a}-vis EC2.

EC1 implies that the ``why'' support sequence for a statement of mathematical knowledge terminates in
finitely many steps at some (set of) reason(s), $X$, which is (are) the most fundamental \emph{when viewed}
\emph{inside the theory at hand}. $X$ is a set of one or more axioms. The restriction accomplished by EC2
is made manifest in practice by choosing $X$ \emph{not} to be an ultimate and only asking questions of the
objects involved which do not require the intake of ultimates. Thus, for the philosophical schools
considered below, mathematical inquiry halts its inspection before reaching the epistemic level of
ultimates. As ultimates are not part of mathematics, mathematics remains mum on any issues 
regarding circularity in the authoritative support (cf. Section \ref{entertaincirc}).
Thus, any circular reasoning arising in mathematics is indeed rightly seen as
problematic.


We will demonstrate for various philosophical schools how EC1 and EC2 are exhibited.\footnote{We 
give only very brief introductions to the schools of thought. For more, see
\cite{KI}; \cite{SW}, especially chapters 17--20; \cite{SH}; \cite{KE}, especially chapter III;
\cite{GG}, Part 5.} We will give special attention to ZFC set theory and the ``Big Three.''\footnote{Term
courtesy of \cite{SH}, Part III.} Our goal is not an exhaustive study of the many philosophies of
mathematics, but a brief introduction to show how each one considered expresses the ECs. In particular,
every philosophy discussed below shares the same question for which the examination would violate EC2,
namely, ``once the foundational elements of EC1 are chosen, \emph{why} are they the correct ones?'' One
could defend the choice of ``these axioms instead of those,'' but the point is that such a defense is a
philosophical argument rather than a mathematical one. Additional comments on EC2 will be mentioned below
as appropriate. Perhaps our discussion will also give the reader a sense for humanity's intellectual
limitations.\footnote{Compare the discussion in \cite[pp.~103-1]{PO}.}

\subsection{ZFC Set Theory: Object Language, Metalanguage, and Russell's Paradox}
The default approach to mainstream mathematics today is termed ZFC set theory, the acronym referring to a
list of 8-9 (depending on the specific formulation) Zermelo-Fraenkel axioms, along with the axiom of
choice.\footnote{We leave an enumeration of the axioms to the references as it is technical and
inessential to the current article. For details, see \cite{SU}; \cite{ST1}; \cite{KI}, especially chapters
VIII--XII; \cite{GG}, especially Entry 5.3; \cite{KE}, especially chapter III; \cite{TA}. For a less
technical introduction, see \cite{PO}, App. E.} The genesis of set theory as its own subfield is often
credited to Georg Cantor and his pioneering work on infinite series, work in which complications were
quickly discovered. \cite[pp.~1-2]{SU} Cantor's work made implicit use of the axiom of abstraction which
asserted the existence of a set whose members are determined by sharing any given property. A 1901 letter
from Bertrand Russell to Gottlob Frege brought to light a fatal flaw with the use of this axiom. To see
the flaw, form the set $S$ whose members share the property that the members are not members of themselves.
Is $S$ a member of itself? If so, then by its own defining property, $S$ is \emph{not} a member of itself.
And if not, then, again by the defining property, $S$ is indeed a member of itself. In both cases, $S$ is
simultaneously a member of itself and not a member of itself. This contradiction was termed ``Russell's
Paradox'' and forced mathematicians to concede the logical flaw in the seemingly innocuous, intuitive
realization of set theory afforded by the axiom of abstraction. \cite[pp.~5-8]{SU}

Russell himself gave an explanation for this state of affairs in 1905.  The problem lies with the use of
\emph{impredicative} definitions, a term first used by Poincaré in 1906, ``wherein an object is
defined (or described) in terms of a class of objects which contains the object being defined. Such
definitions are illegitimate'' since they are circular; hence the related name, Russell's ``vicious
circle principle.'' \cite[pp.~204-7]{KI}, \cite[p.~116]{SH}, \cite[p.~42]{KE}\footnote{Circularity in
such definitions is indeed problematic, even in view of Section \ref{entertaincirc}, since propositional
support is our goal.} Notions such as a ``set of all sets'' are impredicative and hence
barred.\footnote{Some mathematicians and philosophers have defended the use of impredicative definitions.
This debate is beyond the scope of this article. See \cite[pp.~128, 180]{SH}, \cite[pp.~42-5]{KE}. The
Burali-Forti paradox regarding the set of all ordinals is another closely related concept too technical for our coverage.
See \cite[pp.~8, 133]{SU}.} A subtle but powerful modification, introduced by Ernst Zermelo in 1908, led to a
resolution of Russell's Paradox and was included in ZFC as the ``axiom schema of separation.'' It
required that a given set first be known to exist; then, members of that set with a given property can be
collected in a subset. Other paradoxes (or contradictions) were found in early set theory as well, the
details of which we leave to other works. \cite[pp.~5-12]{SU}, \cite[pp.~204-7]{KI}, ch. IX.

What are ZFC's intellectual foundations (EC1)? The axioms and formally undefined terms such as ``set'' are
taken as given and unquestioned.\footnote{Similar comments apply equally well to Russellian type theory,
von Neumann-Bernays-Gödel set theory, and other, alternative set-theoretic foundations.}

How does ZFC exhibit EC2? Through an \emph{object language/metalanguage distinction}.\footnote{Some
authors use the term ``syntax language'' instead of metalanguage. \cite[p.~63]{KE}} Such a distinction is
what allowed the paradoxes to be resolved. The \emph{object language} is the language in which we discuss
the mathematical entities themselves, whether they be sets or numbers or anything else, and the
\emph{metalanguage} is the ``language in which we talk about the object language.'' \cite[p.~11]{SU} The
metalanguage includes the object language, but not vice versa. Indeed, this relationship is what allows
us to leave the notion of set formally undefined in the object language and yet still talk about and
understand it in the metalanguage.  Whereas axioms or definitions are \emph{treated} as ultimates (here, synonymous with ``foundations'') when
doing mathematics, they are only ultimate \emph{with respect to the object language}, not the
metalanguage.\footnote{Kline offers an analogy at \cite[p.~250]{KI}. ``If one wished to study the
effectiveness or comprehensiveness of the Japanese language, to do so in Japanese would handicap the
analysis because it might be subject to the limitations of Japanese.  However, if English is an effective
language, one might use English to study Japanese.''} This distinction makes it possible to ``get
outside'' of set theory (and thus mathematics) to discuss mathematical objects/notions without
\emph{doing} mathematics (and, hence, without being unwittingly restricted by its methods or
limitations). Indeed, much of the history of mathematics since the late 19\textsuperscript{th} century
demonstrates the utility of making this distinction. As Suppes summarizes,

\vspace*{-0.2in}

\begin{quote}
Every one of these paradoxes arises from having available in the language expressions for referring to
other expressions in the language.  Any language with such unlimited means of expression is perforce
inconsistent. Consequently it is important to distinguish between the object language -- here the
language in which we talk about sets -- and the metalanguage, that is, the language in which we talk
about the object language\ldots In other words, we avoid these paradoxes by severely restricting the
richness of our language. \cite[p.~11]{SU}
\end{quote}

\vspace{-2mm}
Let us illustrate this with our $P, Z$, and $X$ notation from above. If $P$ is a statement of knowledge in
mathematics, a complete support sequence of $P$ (a proof of $P$) could be written like ``$P$ because $Q$ because
$R$ because \ldots\ because $X$ because $Y$ because $Z$.'' The supporting reasons up to $X$ would be considered
part of mathematics, and at that point, the mathematical proof (the use of the object language) would
stop.\footnote{Formally, $X$ would be one or more axioms, though in practice most mathematical proofs stop
short of this.} $Y$ and $Z$ would be statements in rational human endeavors which fall outside of mathematics
proper, and thus would be expressions in the metalanguage which are outside of the object language.

\vspace{-1mm}
\subsection{Logicism}

Once the set theory paradoxes were discovered, the question of the consistency (freedom from
contradiction) of mathematics was brought front and center, with the logistic school being the first
serious attempt to address it. \cite[p.~216]{KI} The school credits its founding to Frege in the late
1800s (though its main idea ``can be traced back to Leibniz'' \cite[p.~217]{KI}), and it saw much
development in the early 1900s from mathematicians of renown like Peano, Russell, and Whitehead. The
magnum opus of the school was Russell and Whitehead's \emph{Principia Mathematica}, or \emph{Principia}
for short (which, in particular, included the famous theory of types by which Russell and Whitehead
avoided falling victim to Russell's Paradox). \cite[p.~237]{ST1} ``The purpose of logicism was to show
that classical mathematics is part of logic.'' \cite[p.~122]{SW} This would transfer any foundational
questions about the nature of mathematics into the realm of logic, the consistency of which seemed firmly
established.

\emph{Principia} ``may be considered as a formal set theory'' analogous to ZFC. \cite[p.~122]{SW} Thus
the primary goal became to show that each of the axioms of that set theory, from which flow all of
classical mathematics, are themselves part of (first order) logic. This venture failed since ``at least
two of these axioms, namely, the axiom of infinity and the axiom of choice, cannot possibly be considered
as logical propositions.'' \cite[p.~123]{SW}\footnote{In Russell's formulation, it was the axiom of
reducibility which was contentious.  ``The axiom of reducibility for propositions say {[}sic{]} that any
proposition of higher type is equivalent to one of first order.'' \cite[p.~222]{KI} Russell and
Whitehead admitted that this ``axiom has a purely pragmatic justification: it leads to the desired
results, and to no others.'' \cite[p.~45]{KE}} Mathematicians accept the axiom of infinity, the statement
that infinite sets exist, ``in virtue of its content {[}i.e. from presupposed familiarity with infinite
sets{]} and not in virtue of its syntactical form,'' the latter of which would be necessary to consider
it as a logical proposition. \cite[p.~124]{SW}

Logicism exhibits EC1 by accepting first-order logic as mathematics' starting point. Moreover,
``Frege \emph{proved} that each natural number exists, but his proof is impredicative, violating the type
restrictions {[}of Russell's type theory{]}. Russell had to \emph{assume} the existence of enough
individuals {[}type 0 objects which are neither classes nor sets{]} for each natural number
to exist.'' \cite[p.~119]{SH}

Though logicism is considered to have failed in its primary goal, it remains highly influential in modern
mathematics. ``In fact, it was logicism which started mathematical logic in a serious way.''
\cite[p.~124]{SW}\footnote{For a treatment of mathematical logic, see \cite{TA}. More on the logistic
school can be found in \cite{KI}, chapter X; \cite{KE}, chapter III; \cite{SH}, chapter 5.}

\subsection{Intuitionism}
Lying in stark contrast to logicism is intuitionism, which sees mathematics, even classical mathematics,
very differently. Here, ``mathematics should be defined as a mental activity and not as a set of theorems
(as was done\ldots{[}in{]} logicism).'' \cite[p.~126]{SW} ``Intuition determines the soundness and
acceptability of ideas, not experience or logic.'' \cite[p.~235]{KI} In fact, intuitionism even questions
whether the principles of logic have absolute validity. \cite[p.~46]{KE} Moreover, to be considered
legitimate to intuitionists, the mental activity must be \emph{constructive}.  Objects can only be
considered to exist if one has engaged in every mental step to formulate \emph{completely} that object,
beginning with, for example, the ``primordial intuition'' for the number 1. \cite[p.~125]{SW},
\cite[p.~49]{KE} A precise definition of ``constructive'' is debatable, but, for example, ``Kleene saw
that he could translate it as `compute,' which he \emph{had} precisely defined.'' \cite[p.~iv]{KE}

Intuitionism traces many of its ideas back to Kantian philosophy. \cite[Ch. 7, Sec. 2]{SH} Additionally,
``just as logicism is related to {[}philosophical{]} realism, intuitionism is related to the philosophy
called `conceptualism.' '' \cite[p.~128]{SW} Mathematicians such as Descartes and Pascal had advocated
for a leading role of human intuition, and the mantle was taken up by such luminaries as Kronecker,
Borel, Lebesgue, Poincaré, Baire, and, especially, L. E. J. Brouwer. \cite[pp.~230-1, 234]{KI}

Many mathematicians consider intuitionism unnecessarily restrictive on the tools available to them. For
example, the law of the excluded middle (the principle that either a statement or its negation, but not both,
is true) is considered invalid mathematics by many intuitionists. \cite[p.~126]{SW},
\cite[p.~237]{KI}, \cite[pp.~173-4]{SH}, \cite[pp.~47-8]{KE} The constructive restriction likewise
eliminates the validity of many existence proofs in mathematics. \cite[pp.~238-9]{KI} Further,
constructive proofs are often so much lengthier and more complex that opponents see the burden as too
great. Even more troubling, there are theorems which are true in intuitionistic mathematics but false in
classical mathematics.\footnote{For example, in intuitionistic mathematics, the statement that ``every
real-valued function which is defined for all real numbers is continuous'' is true. \cite[p.~127]{SW}}
This does not bother intuitionists: while the logicists sought to ground all of classical mathematics,
intuitionists do not consider all of classical mathematics legitimate, and the illegitimate parts need no
grounding.  Instead, intuitionists have the goal of giving ``a valid definition of mathematics and then
`wait and see' what mathematics comes out of it,'' fully expecting that some of classical mathematics
will never be realized this way. \cite[p.~126]{SW} ``Brouwer recognized that intuitionistic mathematics
is not a mere restriction of classical mathematics, but is incompatible with it.'' \cite[p.~184]{SH} More
generally, ``as a substitute for classical mathematics {[}intuitionism{]} has turned out to be less
powerful, and in many ways more complicated to develop.'' \cite[pp.~52-3]{KE} For reasons such as these,
this is a philosophical school that the majority of mathematicians have repudiated. And, as with
virtually all philosophical schools, disagreements within intuitionism exist as well. \cite[p.~244]{KI}

``According to intuitionistic philosophy, all human beings have a primordial intuition for the natural
numbers within them\ldots we have an immediate certainty as to what is meant by the number 1\ldots''
\cite[p.~125]{SW} and thus EC1 is manifest in intuitionism. As Leopold Kronecker famously quipped, ``God
made the integers; all the rest is the work of man'' (quoted in \cite[p.~232]{KI}). Brouwer took
mathematical notions like addition and mathematical induction to be intuitively clear, hence
foundational, also. \cite[p.~235]{KI} Even so, intuitionism places mathematics at a higher epistemic
level (i.e. closer to ultimate) than do most other philosophies of mathematics. Since mathematics is
``identical with the exact part of our thinking {[}since it is a mental activity{]}\ldots no science, in
particular philosophy or logic, can be a presupposition for mathematics. It would be [fallaciously] circular to apply
any philosophical or logical principles as means of proof, since mathematical conceptions are already
presupposed in the formation of such principles.\footnote{Intuitionists might therefore revise EC2,
though not dispose of it.} {[}Mathematics has{]} no other source than an intuition,
which\ldots{[}is{]} immediately clear'' (see \cite[p.~51]{KE}, quoting a work from Heyting).

Intuitionism's influence upon the understanding of mathematical existence, as well as on the
related school of constructivism, remain significant today. \cite[p.~244]{KI}, \cite[p.~184]{SH} Also,
``the possibility of an intuitionistic reconstruction of classical mathematics\ldots is not to be ruled
out.'' \cite[p.~53]{KE}

\subsection{Formalism}
As a critic of both the logistic and intuitionist schools, with special vitriol for intuitionism,
\cite[p.~246]{KI} David Hilbert is the primary figure in the formalist school. (That said, many of the
basic notions existed before his ``Hilbert program'' was set in motion around 1910.) The objective of
formalism is, as the name suggests, to formalize a given axiomatization using methods which Hilbert
called ``finitary,'' a somewhat subtle notion which we need not attempt to define here. \cite[pp.~128,
131]{SW}, \cite[p.~250]{KI} The formalization process involves creating a first order language for the
axiomatized theory by enumerating a vocabulary: symbols for items like quantifiers, the equality
relation, connectives (``not,'' ``and,'' ``or,'' etc.), variables, and the undefined terms particular to
the theory. While this process envelops the theory with rigidity, perhaps constraining free-flowing
creativity, it has the advantage of allowing one to study the theory itself via the tools of mathematical
inquiry (an undertaking appropriately termed ``metamathematics'').\footnote{For a thorough discussion of
metamathematics, see \cite{KE}.} Hilbert's purpose was to prove \emph{mathematically} that mathematics
itself was consistent, \cite[pp.~128-30]{SW} and even further, to prove the completeness of the
mathematical enterprise; that is, that ``the axioms\ldots are adequate to establish the correctness or
falsity of any meaningful assertion that involves the concepts'' of mathematics. \cite[p.~158]{KI} For
Hilbert, these twin ventures would establish ``the certitude of mathematical methods'' (quoted
in \cite[p.~158]{SH}).

Philosophically, logicism relates to realism, and intuitionism relates to conceptualism, with the
philosophical leaning of formalism being nominalism. \cite[p.~131]{SW} Nominalism's focus is on the
``names'' of the objects it treats, as its Latin root word suggests. For a formalist, mathematical
objects need not have any existence beyond the name (term, symbol) given to them within the theory. One
consequence is that ``\emph{anything at all} can play the role of the undefined primitives\ldots so long
as the axioms are satisfied.'' \cite[p.~151]{SH}\footnote{In \cite{KI}, Kline includes a similar
discussion of undefined terms in mathematics in the 19\textsuperscript{th} century. The axioms in a
deductive system give an ``implicit definition'' of the terms by telling us how they can be used and what
one knows about them. Thus, the undefined terms in a deductive system ``can be interpreted to be anything
that satisfies the axioms,'' introducing ``a new level of abstraction'' in mathematics; see
\cite[pp.~191ff.]{KI}.} No intuitive meaning is needed. Hilbert, alluding to Euclid's inadequate
definitions,\footnote{For example, Euclid defined a ``point'' as ``that which has no parts.''
\cite[p.~156]{SH} For Hilbert, such vagueness in definition does more harm than good. In such a
situation, ``everything gets lost and becomes vague and tangled and degenerates into a game of hide and
seek.'' Quoted in \cite[p.~156]{SH}. Regarding Euclid, it was also observed (before Hilbert) that in his
\emph{Elements}, ``flaws had arisen because Euclid has been misled by his visual imagery\ldots{[}H{]}e
had inadvertently assumed certain properties of {[}points and lines{]}, without stating them as axioms.''
\cite[p.~243]{ST1}} noted that the ``meaning'' of the undefined terms is given only in ``the
\emph{relations} of points, lines, and planes to each other -- via the axiomatization. All we can provide
is an implicit definition of the terminology'' by specifying their properties through the stated axioms.
\cite[p.~156]{SH}\footnote{As Hilbert put it: ``I have become convinced that the more subtle parts of
mathematics\ldots can be treated with certainty only in this way; otherwise one is only going around in a
circle.'' Quoted in \cite[p.~157]{SH}. Regarding the development of logic and arithmetic, he advocated
for a ``partly simultaneous development,'' essentially a set of co-foundational principles, in order to
avoid ``a vicious circle and\ldots paradoxes.'' Quoted in \cite[p.~111]{SW}.} The formal names can then
be manipulated according to the rules set forth in the first order language, not unlike the way one plays
the game of chess by moving the different pieces according to the game's rules. \cite[p.~131]{SW} This
link is strong enough to have birthed a subschool referred to as ``game formalism.'' \cite[p.~144]{SH} It
is interesting to note that Hilbert developed these notions specifically to avoid circularity in
mathematical reasoning, and for that he required that foundational axioms and definitions be asserted (in
the axiomatization). This suggests
an implicit agreement with the current article's thesis that the rational alternative to foundationalism is
circularity.\footnote{It is even referred to as ``the circular definition problem.'' \cite[p.~122]{BH}}

The failure of formalism to prove the consistency and completeness of mathematics is now famous. Kurt
Gödel's Incompleteness Theorems of the 1930s comprised the fatal blow. \cite[p.~130]{SW} By 1930, some
corners of mathematics had been proven to be consistent and complete (some were even proved by Gödel
himself),\footnote{Notably, first order predicate calculus. See \cite[p.~260]{KI}.} but these results did
not extend to the whole of mathematics. His incompleteness theorems state that, ``if any formal theory
$T$ adequate to embrace the theory of whole numbers is consistent, then $T$ is incomplete. This means
that there is a meaningful\ldots{[}and{]} true statement of number theory which is not provable and so is
undecidable.'' \cite[p.~261]{KI} Putting the Big Three philosophical schools in its crosshairs, ``the
consistency of any mathematical system that is extensive enough to embrace even the arithmetic of whole
numbers cannot be established by the logical principles adopted by the several foundational schools.''
\cite[p.~261]{KI} In other words, ``no axiom system for mathematics as we know it is powerful enough to
lead to a proof of its own consistency.'' \cite[p.~134]{SW}, \cite[p.~167]{SH}

The implications to the mathematical enterprise were far-reaching and led to further unsettling results:
a mathematical proof might have pristine use of logic internally, but there was now no guarantee that
mathematics more broadly was logically meaningful at all. \cite[p.~244]{ST1} It is reported that Hilbert
became furious when he learned of Gödel's work, but that fury has given way in the mathematical community
either to unease or, more commonly, indifference toward the philosophy of the foundations of mathematics.
\cite[p.~245]{ST1}, \cite[p.~132]{SW} A fuller account of these issues (crises!) will be left to the
references.\footnote{See \cite{KI}, chapter XII; \cite{SH}, sections 6.4, 8.1; \cite{KE}, section 42. The
notion of (un)certainty in mathematics is a broader theme of many of these works; cf. \cite[pp.~212-6]{BH}.}

What about the epistemic characteristics' relationship to formalism? On EC1, formalism begins with an
axiomatized theory (a set of axioms and undefined terms), so the items in such a theory are deemed to be
supplied beforehand.

Regarding EC2, ``one has to talk about the {[}first order{]} language $L$ as one object, and while doing
this, one is not talking within that safe language $L$ itself, \ldots{} {[}but in{]} ordinary, everyday
language\ldots{[}In so doing,{]} there is of course every danger that contradictions, in fact, any kind
of error, may slip in.'' \cite[pp.~130-1]{SW}\footnote{Notice the object language/metalanguage
distinction.} Hilbert's solution was the aforementioned insistence on finitary reasoning, a restriction
which demonstrates EC2. EC2 is displayed perhaps even more clearly by Gödel's rebuke of the formalist
school.  Whereas Gödel proved that a mathematical system cannot show \emph{its own} consistency, at least
some of the undecidable statements in a formal system can be shown to be true by \emph{informal}
arguments. \cite[p.~263]{KI} As Gödel put it, ``it is necessary to go beyond the framework of what is, in
Hilbert's sense, finitary mathematics if one wants to prove consistency of classical
mathematics.''\footnote{Quoted in \cite[p.~167]{SH}.} And again, ``the methods which we must trust in the
proof {[}of consistency{]} must include some which lie outside the collection of the methods formalized
in the system.'' \cite[p.~211]{KE} Natural, interesting, and even crucial questions about mathematics
must use tools which transcend the discipline. Since such tools evidently exist, Gödel's theorems
substantiate the claim that mathematical knowledge is not ultimate knowledge.

Despite the formalist's failure in their initial goal, the school has had immense, ongoing impact on the
way mathematics is done today. Fields such as metamathematics or proof theory,\footnote{These terms are
introduced in \cite[p.~55]{KE}.} as well as ``modern mathematical logic and its various offshoots, such
as model theory, recursive function theory, etc.'' \cite[p.~130]{SW} owe an immeasurable debt to
formalism. The emphasis on stripping objects of their ``unnecessary'' features (whatever is not pertinent
to the logical structure), and the recognition that physical intuition can be misleading, help explain
why ``so much of the abstraction of early 20\textsuperscript{th} century mathematics stemmed from
Hilbert's viewpoint.'' \cite[p.~243]{ST1}, \cite[pp.~54-5]{KE}

\subsection{Other Philosophies of Mathematics}
For each of the philosophical schools given above, there remain adherents of (modified versions of)
them,\footnote{Shapiro provides examples throughout \cite{SH}.} but they don't exhaust mathematicians'
and philosophers' present-day views. We will make passing mention of some other selected philosophies of
mathematics, some of which overlap the above (and each other). As described by
theologian and mathematician Vern Poythress, ``Platonism says that numbers and mathematics belong to a
realm of abstract ideas, a realm that exists before mathematicians begin to study it.'' \cite[p.~153]{PO}
An alternative view, empiricism, seeks to found mathematics in the experiences of the senses.
\cite[pp.~154-5]{PO} Predicativism ``accepts the natural numbers as given\ldots either by our intuition
or by a Platonic realm or by both,'' and then seeks to rebuild as much of mathematics as possible while
avoiding impredicative definitions. \cite[p.~159]{PO}\footnote{For a thorough, favorable treatment of
predicativism, see \cite{SR}. Information available at \url{https://philpapers.org/rec/STOADO-2}.
Accessed July 12, 2022.} These three views each begin mathematics proper with (at least) the natural
numbers, and thus all exhibit EC1.\footnote{Their ontological claims about the nature of
mathematical objects differ, but that is immaterial to our discussion.} William van Orman Quine's philosophical
naturalism advocates for science as the guide for what constitutes acceptable mathematics,
displaying EC1 by assuming results of scientific inquiry before mathematics can
commence. \cite[p.~160]{PO}

Shapiro's treatment of some contemporary philosophies of mathematics tells a similar story.
\cite[Part~IV]{SH} The core difference between the philosophers and their associated schools (e.g.
realism vs.  anti-realism) is the ontology of mathematical objects. But in each case, \emph{something} is
provided to mathematics (frequently, the natural numbers) before it can begin. EC1 thus appears in each.
As two examples: in structuralism, the structures themselves which mathematicians study are the givens.
\cite[pp.~258ff.]{SH} In Hartry Field's fictionalism, specifically his realization of Newtonian
gravitational theory, ``points'' and ``regions'' are (among) the givens. \cite[pp.~229ff.]{SH}

It is worth noting that Poythress's \emph{Redeeming Mathematics} proposes (with remarkable succinctness
and in lay terminology) an alternative philosophy of mathematics which could, in some sense, be
considered an exception among the philosophies of mathematics. He grounds the tools of mathematics in the
Christian God, who is ultimate in the Christian worldview. Poythress criticizes the ``reductionisms'' of
the aforementioned philosophies of mathematics\footnote{ ``Each of these {[}philosophies{]} has a
preferred starting point\ldots{[}which{]} becomes the preferred platform for explaining everything else
in mathematics.'' \cite[p.~151]{PO} Note the close connection of this observation with EC1.} and lays out
a more holistic approach which might require refinement of the ECs.\footnote{The comments in
\cite[pp.~104-5]{PO} are prime examples. In Poythress' view, God, and hence ultimate knowledge, is
inextricably linked to mathematics. The pivotal feature is then the internal coherence of the system, not
which ``starting point'' is chosen.}

\subsection{(The Lack of) The Epistemic Characteristics in Theology and Philosophy}\label{lackECtheophil}
If mathematics has been able to employ EC1 and EC2 with great success (as have other academic
disciplines), can fields such as theology, philosophy, religion, or worldview studies emulate it?
These latter fields are of a fundamentally different character to the extent that duplicating the pattern
is not possible. Many, though not all, of their questions are inherently about ultimates (recall we
refer to these as ``ultimate questions''). For the ultimate questions in these fields (or in whatever
field one places ultimate questions), there is no way to exhibit both EC1 and EC2.


Most significant for the present discussion is in regards to EC2.\footnote{To briefly address EC1: while different philosophical or theological schools of thought disagree on the starting points of EC1, each school presumably has some such foundational set. For example, see \cite {GR} for a comparison of starting points for generalists and particularists, and \cite{CH} for a thorough treatment of the deep epistemological issues faced in determining the elements which comprise EC1 (note the assumption that such elements exist).} Limiting the scope of ultimate questions is not possible, for to do so would be to render
them non-ultimate. In this case, an object language/metalanguage distinction is impossible; that is, the
object language also includes the metalanguage, and thus the two are equal (recall that the reverse
inclusion holds for mathematics and any other discipline in which the distinction is
utilized).\footnote{This is illustrated in simple terms in \cite{SW}: ``when the logicists tell us what
they mean by a logical proposition, \ldots{} they use philosophical and not mathematical language. They
have to use philosophical language for that purpose since mathematics simply cannot handle definitions of
so wide a scope.'' \cite[p.~124]{SW} Similar comments for intuitionism can be found in
\cite[pp.~127-8]{SW}.} There is not, and indeed cannot be (for our human agent $H$), a greater context in which to
discuss them. Since these questions arise naturally in theology and philosophy, there is no way to ``get
outside'' of those disciplines when investigating such questions; one cannot attempt to answer ultimate
questions in theology without theologizing. There is nothing in the line of reasoning that falls outside
of the umbrella of, or transcends, the discipline. As Van Til put it, ``{[}w{]}e must go round and round
a thing to see of its dimensions and to know more about it, in general, unless we are larger than that
which we are investigating,'' \cite[p.~24]{V1}, \cite[p.~257]{FR} and in the case of ultimate questions,
we are by definition not ``larger'' than they are.  In this sense, theology and philosophy are
intellectually all-encompassing. Evidently it is this observation, in part, which motivates worldview
scholar Nancy Pearcey to refer to her Christian worldview as ``Total Truth.'' \cite{PE} EC2 cannot be
exhibited for ultimate questions. 

Ultimates thus have a ``meta'' relationship to other statements of knowledge.  They don't simply lead to
or propositionally support other statements; rather, they form the framework in which an intellectual discipline is done.
Intelligibility itself, a precondition of any intellectual pursuit, is determined by the
ultimates.\footnote{In theology in particular, ultimates are often referred to as presuppositions.
On this topic, Frame says, ``The `pre' {[}in presupposition{]} should be understood mainly as an indicator
of eminence\ldots not temporal priority.'' \cite[p.~xxxii]{FR}}
Since it is fields such as philosophy and theology which handle ultimate questions, the ECs cannot both
be always exhibited therein.\footnote{Note that some of the questions which are
asked in theology, philosophy, and related disciplines \emph{can} be handled in such a way that both EC1 and EC2 can be exhibited.  The main point is that, unlike mathematics, \emph{not all} questions in
theology or philosophy can exhibit \emph{both} epistemic characteristics (cf. the discussion of
EC1 in subsection \ref{Godneedcirc}).} One particular conclusion, in view of (FC), is that if circular support is acceptable anywhere, it would be so in fields such as theology and philosophy by virtue of their inclusion of ultimate questions.\footnote{Recall that it is precisely the exhibition of the ECs which allows mathematics to categorically reject circular reasoning; see Section \ref{TheECs}, and compare the nuance in the meaning of ``circularity'' in Section \ref{entertaincirc}.} 
%

\subsection{Analogies}
As an epilogue to the primary thrust of the article, we will make brief mention of a potpourri of
analogies, both within and without mathematics.

First, the epistemic distinction between ultimates and non-ultimates bears a resemblance to the
distinction in properties between the boundary and interior of certain mathematical objects, such as a
manifold. A planar disk (a circle with all of the interior points) is an example. The interior
has different properties (most notably, it is 2-dimensional) than the boundary (the 1-dimensional edge of
the disk). As ultimates are at the ``boundary'' of human knowledge, it is not surprising that their
epistemological properties differ qualitatively from the ``interior,'' non-ultimate knowledge.

Consider another mathematical analogy: a local/global distinction. A local property is satisfied only at
or within a neighborhood of a particular point/object (for example, any open neighborhood of a point on a
2-sphere smaller than the 2-sphere itself, is isomorphic to the real plane), whereas a global property
in mathematics refers to the set of all points (objects) under consideration (the 2-sphere, taken as an
entire object, is \emph{not} isomorphic to the real plane). The realm of epistemological influence of a non-ultimate
proposition can be thought of as local, whereas an ultimate has global influence.

Many of the above concepts are not unique to mathematics, so it can be illuminating to see an
``everyday'' example. One appears in dictionaries: a word cannot be defined without the use of
other words.  Since there are only finitely many words, this means that, eventually, some word has to be
left either undefined (deemed so foundational as to be left out of the dictionary), or there is a sequence
of words which are used in
each other's definition (which is circularity; examples can indeed be found in modern dictionaries).
(FC) applies to definitions as well as propositions.

The object language/metalanguage distinction can also be seen outside of mathematics. Consider a
theatrical script: it employs stage directions (metalanguage) in addition to the lines of the show
(object language).  Similarly, some government or legal documentation says, ``This page
intentionally left blank.'' Those words themselves cause the page not to be blank, but they are
considered part of the document's metalanguage rather than the document's content (object language).

\section{Follow-up Questions}\label{followup}
\subsection{Does this article itself use foundational or circular reasoning?}\label{circinarticle}
This article uses foundationalism in its propositional progression. 
We assumed several notions at the outset, such as the use
of logic and rationality in our truth-seeking. These were taken as ultimate (without further defense)
\emph{for this article}. This particular assumption illustrates the ``meta'' relationship that ultimates
have with respect to other statements: had we not assumed anything regarding logic and rationality, by what standards would
we have proceeded to answer any questions or make any claims? We could not reason about reason without
using reason! The rationality assumption established a precondition for the remainder of the article.


\subsection{Could different individuals choose different ultimates, resulting in diverse but
self-consistent worldviews that are mutually exclusive?}
In mathematics, this circumstance (a panoply of axiomatic theories) was essentially endorsed by David
Hilbert. ``Literally, Hilbert claimed that if a collection of axioms is consistent, then they are true
and the things the axioms speak of exist.'' \cite[p.~156]{SH} Axioms in mathematics can be arbitrarily
chosen so long as they are consistent - a version of coherence epistemology in mathematics.
In fact, it is this line of reasoning which allows for the
parallel\footnote{No pun intended.} existence of multiple geometries (Euclidean and non-Euclidean) by
asserting different axioms.  Mathematicians are comfortable with this state of affairs, so can we apply
the same reasoning to worldviews, inclusive of ultimates?

A priori an affirmative answer to this question seems possible, and by our assumptions it seems to be so
\emph{in this article.} However, let us momentarily suspend our moratorium on evaluating theories of warrant.
It is the author's belief that, in practice, all truth-conscious
individuals operate, implicitly or explicitly, under some shared presuppositions\footnote{Presuppositions
in this context are ultimates, but as the former term is used more commonly when discussing worldviews,
we will adopt it for this follow-up question.} which eliminate some worldviews from rational contention; in
particular, the possibility of global mutual exclusivity is removed.  These presuppositions include
shared ultimate standards (e.g. for rationality and morality) by which everything else is evaluated.
The commonality between persons' presuppositions assuages fears that two individuals could
have perfectly consistent and veracious worldview theories and yet have nothing on which they agree; in reality, the
two don't actually inhabit disjoint systems.

For example, any theological or philosophical system that includes the principle of noncontradiction, as
is assumed here, requires internal consistency.
A necessary condition for all other ultimate criteria would be that they not refute
anything already in the system of ultimates. This condition alone can indeed rule out many worldviews which
might otherwise seem to be viable a priori. Details and a more thorough treatment of how to choose the ``right''
presuppositions, or what it might mean to be the right presuppositions, would be the topic of another
work.

\subsection{Does God need foundationalism, circularity, both, or neither?}\label{Godneedcirc}
One would first have to define what is meant by ``God.'' Specifically in Christianity and the other
Abrahamic faiths, the conception of God is traditionally as an \emph{ultimate} God (in particular, He is
or determines the ultimate standard for rationality and morality). In other words, by definition, there is no higher
standard by which God can be judged. 
His word can be taken as foundational, and His authority (and anything in His nature) as self-verifying; 
the notion of circularity in \emph{authoritative} support seems to be consistent with God's nature in these faiths. This
concept is part of what theologians refer to as divine aseity,\footnote{``Something `exists necessarily'
if it cannot fail to exist. Something `exists contingently' if it can fail to exist. Whether or not
something contingent exists depends on factors outside itself. Therefore, `contingent existence' (such as
ours) is `dependent existence.' Necessary existence is aseity or self-existence.'' \cite[p.~116n36]{FR}
For more on divine aseity, see \cite[pp.~265-78]{FR}. For more on the relationship between divine aseity
and mathematics, see \cite[pp.~63ff.]{PO}} a concept hinted at in Hebrews 6:13;\footnote{``For when God
made a promise to Abraham, since he had no one greater by whom to swear, he swore by himself\ldots''
\cite{ESV}} God has no standard to appeal to above Himself! Thus, anything that is known to be part of
God's nature need not (indeed cannot) be further justified.\footnote{The sticking point is knowing which
characteristics are part of His nature. There is not agreement there, even among devoted practitioners
of a particular faith tradition. But the point is that if the qualities of His nature are assumed to be known, then
those qualities ground completely the authoritative justification. A common objection is that grounding a standard like rationality or
morality in God's nature is ``kicking the can down the road,'' suggesting that there is another step of
reasoning beyond God's nature which requires explanation. Such an objection violates the conception of
God as being Himself ultimate (by definition), or it seems to make the implicit assumption that finite
termination of the support sequence is impossible. Either is problematic.} 
We must be careful, though,
because it is quite possible that the assumptions of the current article do not all necessarily apply to
God's ``divine knowledge.'' In particular, the finiteness and time-bound nature of humanity are evidently not limitations for
God. Thus, we do not claim that the argumentation of the current article applies to divine
knowledge, nor is a claim made excluding the possibility of unimagined spiritual or
supernatural realities about which humans have never received revelation nor have
other means of gaining knowledge.

One final comment on the relationship between Christian faith and human knowledge is in order. The Bible,
especially the wisdom literature (such as the book of Proverbs), affirms that human knowledge,
mathematics included, finds its final foundation in (reverence for) God. Job 28:28; Psalm 111:10;
Proverbs 1:7, 2:6, 3:5-6, 8:22, 9:10; and Colossians 2:2-3 are examples.\footnote{They read as follows. 
``And he said to man,
`Behold, the fear of the Lord, that is wisdom, and to turn away from evil is understanding' '' (Job
28:28). ``The fear of the Lord is the beginning of wisdom; all those who practice it have a good
understanding. His praise endures forever!'' (Psalm 111:10). ``The fear of the Lord is the beginning of
knowledge; fools despise wisdom and instruction'' (Proverbs 1:7).  ``For the Lord gives wisdom; from his
mouth come knowledge and understanding'' (Proverbs 2:6). ``Trust in the Lord with all your heart, and do
not lean on your own understanding. In all your ways acknowledge him, and he will make straight your
paths'' (Proverbs 3:5-6). ``{[}Wisdom says{]} `The Lord possessed me at the beginning of his work, the
first of his acts of old' '' (Proverbs 8:22). ``The fear of the Lord is the beginning of wisdom, and the
knowledge of the Holy One is insight'' (Proverbs 9:10). ``\ldots that their hearts may be encouraged,
being knit together in love, to reach all the riches of full assurance of understanding and the knowledge
of God's mystery, which is Christ, in whom are hidden all the treasures of wisdom and knowledge''
(Colossians 2:2-3). Scripture quotations are from \cite{ESV}.} At risk of over-enthusiastic hermeneutics,
verses like these seem to suggest that, in a Christian worldview, EC1 will eventually be exhibited by all
human knowledge, with the starting point being God's own self-authenticating nature. In the author's
view, the philosophical skepticism of Agrippa and the Pyrrhonists (and many philosophers of our day) is
indeed correct that human knowledge would be hopeless without a ``beginning;'' and moreover, that that
beginning must be external to humanity in order to be ultimately trustworthy. By virtue of the doctrines of
historic Christian belief, the nature of God is sufficient for such a beginning.

The author would like to thank the anonymous reviewers for comments which improved the article.

\setbiblabelwidth{100}


\begin{thebibliography}{99}

\bibitem{AL} Albert, Hans. \emph{Treatise on Critical Reason}, (Princeton, NJ: Princeton University
Press), 1985. Translated by Mary Varney Rorty. Chap. 1 Sect.  2.

\bibitem{B1} Bahnsen, Greg L. \emph{Always Ready: Directions for Defending the Faith}, (Nacogdoches, TX:
Covenant Media Press), 1996.

\bibitem{B2} Bahnsen, Greg L. \emph{Van Til's Apologetic: Readings and Analysis}, (Phillipsburg, NJ: P\&R
Publishing), 1998.

\bibitem{BH} Bradley, James and Howell, Russell. \emph{Mathematics Through the Eyes of Faith}, (New York:
Harper Collins), 2011.

\bibitem{CZ} Chartrand, Gary and Zhang, Ping. \emph{A First Course in Graph Theory}, (Mineola, NY: Dover Publications,
Inc.), 2012.

\bibitem{CH} Chisholm, Roderick M. \emph{The Problem of the Criterion}, (Milwaukee, WI: Marquette University Press), 1973.

\bibitem{CL} Clark, Kelly James. \emph{Return to Reason}, (Grand Rapids, MI: Eerdmans Pub. Co.), 1990.

\bibitem{EM} Empiricus, Sextus. \emph{Outlines of Scepticism}, (Cambridge: Cambridge University Press),
2000.

\bibitem{EN} Engel, S. Morris. \emph{With Good Reason: An Introduction to Informal Fallacies,
5\textsuperscript{th} ed.}, (New York: St. Martin's Press), 1994.

\bibitem{EG} Engel, Jr., Mylan. \emph{Positism: The Unexplored Solution to the Epistemic Regress
Problem}, Metaphilosophy, Vol. 45 No. 2, April 2014.

\bibitem{FR} Frame, John M. \emph{Apologetics: A Justification of Christian Belief}, (Phillipsburg, NJ:
Presbyterian and Reformed Publishing Co.), 2015.

\bibitem{GG} Grattan-Guinness (ed.). \emph{Companion Encyclopedia of the History and Philosophy of the
Mathematical Sciences}, (London: Routledge), 1994.

\bibitem{GR} Greco, John. \emph{Methodology in epistemology: particularism and generalism},
available at \url{https://www.rep.routledge.com/articles/overview/epistemology/v-3/sections/methodology-in-epistemology-particularism-and-generalism}. Accessed March 21, 2023.

\bibitem{HA} Haack, Susan. \emph{Evidence and Inquiry: Towards Reconstruction in Epistemology}, (Oxford:
Blackwell Publishers), 1995.

\bibitem{ESV} The Holy Bible, English Standard Version (ESV): Containing the Old and New Testaments (Wheaton, IL: Good
News Publishers), 2001.

\bibitem{HB} Howell, Russell W. and Bradley, W. James (eds.). \emph{Mathematics in a Postmodern Age: A Christian
Perspective}, (Grand Rapids, MI: Wm. B. Eerdmans Publishing Co.), 2001.

\bibitem{KE} Kleene, Stephen Cole. \emph{Introduction to Metamathematics}, (NY: Ishi Press
Int'l), 2009.

\bibitem{KI} Kline, Morris. \emph{Mathematics: The Loss of Certainty}, (NY: Oxford Univ.
Press), 1980.

\bibitem{PE} Pearcey, Nancy. \emph{Total Truth: Liberating Christianity from Its Cultural Captivity},
(Wheaton, IL: Crossway Books), 2005.

\bibitem{PL} Plantinga, Alvin. \emph{Reason and Belief in God}, in Plantinga and Wolterstorff, Nicholas (eds.).
\emph{Faith and Rationality}, (Notre Dame, IN: University of Notre Dame Press), 1983, pp. 16-93.

\bibitem{PO} Poythress, Vern S. \emph{Redeeming Mathematics: A God-Centered Approach}, (Wheaton,
Illinois: Crossway), 2015.

\bibitem{PR} Pratt, Richard L., Jr. \emph{Common Misunderstandings of Van Til's Apologetics}, part 2,
available at \url{http://www.thirdmill.org/newfiles/ric_pratt/TH.Pratt.VanTil.2.html}. Accessed March 21, 2023.

\bibitem{SO} Sosa, Ernest. \emph{The Raft and the Pyramid: Coherence versus Foundations in the Theory
of Knowledge}, (Midwest Studies in Philosophy) Vol. 5 No. 1, 1980, pp. 3-26.

\bibitem{SP} Sproul, R. C. \emph{Scripture Alone: The Evangelical Doctrine}, (Phillipsburg, NJ: P\&R), 2005.

\bibitem{SH} Shapiro, Stewart. \emph{Thinking About Mathematics}, (New York: Oxford University Press),
2000.

\bibitem{ST1} Stewart, Ian. \emph{The Story of Mathematics From Babylonian Numerals to Chaos Theory},
(London: Quercus Publishing), 2008.

\bibitem{SR} Storer, T. \emph{A defence of predicativism as a philosophy of mathematics}, (Doctoral
thesis), 2010.

\bibitem{SU} Suppes, Patrick. \emph{Axiomatic Set Theory,} (New York: Dover Publications, Inc.), 1960.

\bibitem{SW} Swetz, Frank (ed.). \emph{The Search for Certainty: A Journey Through the History of
Mathematics from 1800--2000}, (Mineola, NY: Dover Publications, Inc.), 2012.

\bibitem{TA} Tarski, Alfred. \emph{Introduction to Logic and to the Methodology of the Deductive
Sciences}, (New York: Dover Publications, Inc.), 1995. Reprinting of the original 1946 textbook.

\bibitem{TO} Torres, Joseph E. \emph{Between Scylla and Charybdis: Presuppositionalism, Circular
Reasoning, and the Charge of Fideism Revisited}, included as Appendix D in Frame, 2015.

\bibitem{V1} Van Til, Cornelius. \emph{The Metaphysics of Apologetics}, (unpublished class syllabus),
1932.

\bibitem{V2} Van Til, Cornelius. \emph{A Survey of Christian Epistemology}, (Philadelphia: P\&R), 1969.

\bibitem{W1} Walton, Douglas N. \emph{Begging the Question: Circular Reasoning as a Tactic of
Argumentation}, (New York: Greenwood Press), 1991.

\bibitem{W2} Walton, Douglas N. \emph{Circular Reasoning}, in Blackwell Companion to Epistemology, ed.
Jonathan Dancy and Ernest Sosa (Cambridge, MA: Blackwell Reference), 1992.

\end{thebibliography}
\end{document}